\documentclass[a4paper,11pt]{article}

\usepackage[vcentermath]{youngtab}

\usepackage[dvips]{graphicx}
\usepackage{amsmath, epsfig, cite}
\usepackage{amssymb}
\usepackage{amsfonts}
\usepackage{latexsym}
\usepackage{graphicx}
\usepackage{extarrows}
\usepackage{booktabs}
\usepackage{color}
\usepackage{tikz}
\usetikzlibrary{calc,through,backgrounds}

\makeatletter

\newcommand{\Rmnum}[1]{\expandafter\@slowromancap\romannumeral #1@}
\makeatother

\def\proof{\noindent {\it{Proof.} \hskip 2pt}}

\newtheorem{thm}{Theorem}[section]

\newtheorem{lem}[thm]{Lemma}
\newtheorem{conj}[thm]{Conjecture}

\newcommand{\qed}{{\hfill\rule{4pt}{7pt}}}

\numberwithin{equation}{section}

\newcommand{\latticepath}[4]{ 
  \coordinate (L) at #2;
  \foreach \x/\y/\a in {#4} {
    \coordinate (L1) at ($ (L) + ( #1 * \x , #1 * \y ) $);
    \draw[color=\a, #3] (L) -- (L1);
    \coordinate (L) at (L1);
  }
}

\newcommand{\tikzbox}[8]{ 
  \coordinate (A) at #2;
  \coordinate (B) at ($ (A) + (  #1,0    ) $);
  \coordinate (C) at ($ (A) + (   0,-#1  ) $);
  \coordinate (D) at ($ (A) + (  #1,-#1  ) $);
  \coordinate (E) at ($ (A) + (#1/2,-#1/2) $);
  \draw[fill=#7] (A) rectangle (D);
  \draw[color=#3] (A) -- (B);
  \draw[color=#6] (A) -- (C);
  \draw[color=#4] (B) -- (D);
  \draw[color=#5] (C) -- (D);
  \node at (E) []{#8};
}

\newcommand{\blackbox}[3]{ 
  \tikzbox{#1}{#2}{black}{black}{black}{black}{white}{#3}
}

\newcommand{\boxcollection}[3]{ 
  \coordinate (X) at #2;
  \foreach \x/\y/\object in {#3} {
    \blackbox{#1}{($ (X) + (#1 * \y, - #1 * \x) - (#1, - #1) $)}{\object};
  }
}

\newcommand{\fl}[1] {\lfloor #1 \rfloor}

\newcommand{\fs} {\fl {\frac s 2}}
\newcommand{\ft} {\fl {\frac t 2}}

\makeatletter \@addtoreset{equation}{section} \makeatother

\begin{document}
\rule{0cm}{1cm}
\begin{center}
{\Large\bf Average Size of a Self-conjugate $(s,t)$-Core Partition}
\end{center}
 \vskip 6mm
 \begin{center}
{\small William Y.C. Chen$^1$,
Harry H.Y. Huang$^2$ and Larry X.W. Wang$^3$}

\vskip 4mm
$^{1,2,3}$Center for Combinatorics, LPMC-TJKLC\\
Nankai University,
Tianjin 300071,
P.R. China \\[3mm]

\vskip 4mm

$^1$chen@nankai.edu.cn,   $^2$hhuang@cfc.nankai.edu.cn   \\
$^3$wsw82@nankai.edu.cn,
\end{center}

\begin{abstract} Armstrong, Hanusa and Jones conjectured that if $s,t$ are coprime integers, then the average size of an $(s,t)$-core partition and the average size of a self-conjugate $(s,t)$-core partition are both equal to $\frac{(s+t+1)(s-1)(t-1)}{24}$. Stanley and Zanello showed that  the average size of an $(s,s+1)$-core partition equals $\binom{s+1}{3}/2$.  Based on a bijection of  Ford, Mai and Sze  between self-conjugate $(s,t)$-core partitions and lattice paths in $\lfloor \frac{s}{2} \rfloor\times \lfloor \frac{t}{2}\rfloor$ rectangle, we
obtain the average size of  a self-conjugate $(s,t)$-core partition as conjectured by Armstrong, Hanusa and Jones.

\end{abstract}

\noindent {\bf Keywords:} $(s,t)$-core partition, self-conjugate partition, lattice path

\noindent {\bf AMS Classification:} 05A17, 05A15

%

\section{Introduction}



In this paper, employing a bijection of Ford, Mai and Sze between self-conjugate $(s,t)$-core partitions and lattice paths, we prove a conjecture of Armstrong, Hanusa and Jones on the average size of a self-conjugate $(s,t)$-core partition.


A partition is called a $t$-core partition, or simply a $t$-core,
if its Ferrers diagram contains no cells with hook length $t$.
A partition  is called an $(s,t)$-core partition, or simply an $(s,t)$-core,
if it is simultaneously an $s$-core and a $t$-core. When $gcd(s,t)=r>1$, each $r$-core is an $(s,t)$-core, which means that there are infinitely many $(s,t)$-cores.
When $s$ and $t$ are coprime, Anderson \cite{jaclynAnderson} showed that the number of $(s,t)$-core partitions equals \[\frac{1}{s+t}{{s+t}\choose s}.\]
Under the same condition,
Ford, Mai and Sze \cite{Ford_Mai_Sze} characterized the 
set of hook lengths of diagonal cells in self-conjugate $(s,t)$-core partitions,
and they showed that the number of self-conjugate $(s,t)$-core partitions is
\begin{equation}\label{SCNumber}
{{{\lfloor \frac s 2 \rfloor}+{\lfloor \frac t 2 \rfloor}}\choose{{\lfloor \frac s 2 \rfloor}}}.
\end{equation}

A partition is of size $n$ if it is a partition of $n$.  Aukerman, Kane and Sze \cite{aks} conjectured that the largest size of an $(s,t)$-core partition for $s$ and $t$ are coprime. Olsson and Stanton \cite{OlssonStanton} proved this conjecture and gave the following stronger theorem.
\begin{thm}\label{largest}
If $s$ and $t$ are coprime, then there is a unique largest $(s,t)$-core partition (which happens to be self-conjugate) of size
\begin{eqnarray}\label{largestSize}
\frac{(s^2-1)(t^2-1)}{24}.
\end{eqnarray}
\end{thm}
A short
proof for the conjecture of Aukerman, Kane and Sze was given by Tripathi \cite{Tripathi}. Vandehey \cite{kid} gave the following characterization of the largest $(s,t)$-core partition.

\begin{thm}\label{contain}
There exists a largest $(s,t)$-core partition $\lambda$ under the partial order of containment. That is, for each $(s,t)$-core $\mu$, $\lambda_i\geq \mu_i$ for $1\leq i \leq l(\mu)$.
\end{thm}

It is clear that the largest $(s,t)$-core in the above theorem is unique.
It is
 the $(s,t)$-core of the largest size, and it is also a $(s,t)$-core
 of the longest length.

Recently, Armstrong, Hanusa and Jones \cite{Arm1308} proposed the following conjecture
concerning the average size of an $(s,t)$-core and the average size of a self-conjugate $(s,t)$-core.

\begin{conj}
  {\label{armConj}}
Assume that $s$ and $t$ are coprime. Then the average size of an $(s,t)$-core and
 the average size of a self-conjugate $(s,t)$-core are both equal to
\begin{equation*}    
\frac{(s+t+1)(s-1)(t-1)}{24}.
\end{equation*}
\end{conj}

Stanley and Zanello \cite{StanleyZanello} showed that the conjecture for
the average size of an $(s,t)$-core holds for $(s,s+1)$-cores.
More precisely, they obtained that the average size of an $(s,s+1)$-core equals $\binom{s+1}{3}/2$. In this paper, we prove the Conjecture \ref{armConj} pertaining to
 the average size of a self-conjugate $(s,t)$-core.

%

\section{Proof of the conjecture for self-conjugate $(s,t)$-cores}

In this section, we prove the conjecture of Armstrong, Hanusa and Jones for self-conjugate $(s,t)$-cores.
Let us begin with a quick review of the work on the structure of self-conjugate $(s,t)$-cores. Define $$MD(\lambda):= \{h | h  \text{ is the hook length of a cell on the main diagonal of } \lambda \}.$$
It is easily seen that a self-conjugate partition is uniquely determined by its main diagonal hooks. Ford, Mai and Sze \cite{Ford_Mai_Sze} gave the following characterization of the main diagonal hook length set of a self-conjugate $t$-core $\lambda$.

\begin{thm}\label{MD}
A self-conjugate partition $\lambda$ is a $t$-core if and only if  the following
 conditions hold: \newline
(1) if  $h \in MD(\lambda)$ and $h > 2t$, then $h-2t$ is also in $MD(\lambda)$;  \newline 
(2) if $h,s \in MD(\lambda)$, then $h+s \neq 0 \pmod{2t}$.
\end{thm}


To characterize the main diagonal hook lengths of a self-conjugate $(s,t)$-core, Ford, Mai and Sze \cite{Ford_Mai_Sze} introduced
an integer array $A=(A_{i,j})_{1\leq i \le \fl{s/2},1\leq j \leq \fl{t/2}}$, where
\begin{equation}\label{aij}
A_{i,j}= st - (2j-1)s-(2i-1)t,
\end{equation}
for $1\leq i \le \fl{s/2}$ and $1\leq j \leq \fl{t/2}$.
Let $\mathcal{P}(A)$ be the set of  lattice paths  in $A$ from  the lower-left corner to the upper-right corner. For example,
Figure \ref{PathInArray} gives an array $A$ for $s=8$ and $t=11$, and the solid lines
represent a lattice path in $\mathcal{P}(A)$.
For a lattice path $P$ in $\mathcal{P}(A)$, let $M_A(P)$  denote the set of positive entries $A_{i,j}$ below $P$ and the absolute values of negative entries above $P$.

The following theorem is due to Ford, Mai and Sze \cite{Ford_Mai_Sze}.

\begin{thm}\label{thmCorner}
Assume that $s$ and $t$ are coprime. Let $A$ be the array as given  in  (\ref{aij}).
 Then there is a bijection
 $\Phi$  between the set $\mathcal{P}(A)$ of lattice paths and
 the set of self-conjugate $(s,t)$-core partitions such that
 for $P\in \mathcal{P}(A)$, the set of main diagonal hook lengths of $\Phi(P)$
 is given by $M_A(P)$.
\end{thm}

For example, for the lattice path $P$ in Figure \ref{PathInArray}, $5$ is the only positive entry below $P$, while $-7$ and $-13$ are the negative entries above $P$.  Thus $M_A(P)=\{5,7,13\}$. This gives $\Phi(P)=(7,5,5,3,3,1,1)$, which is an
 $(8,11)$-core partition.

\begin{figure}
  \centering
  \begin{tikzpicture}[scale=1]

    \boxcollection{1}{(0,0)}{
1 / 1 / $ 69 $ ,1 / 2 / $ 53 $ ,1 / 3 / $ 37 $ ,1 / 4 / $ 21 $ ,1 / 5 / $ 5 $ ,
2 / 1 / $ 47 $ ,2 / 2 / $ 31 $ ,2 / 3 / $ 15 $ ,2 / 4 / $ -1 $ ,2 / 5 / $ -17 $ ,
3 / 1 / $ 25 $ ,3 / 2 / $ 9 $ ,3 / 3 / $ -7 $ ,3 / 4 / $ -23 $ ,3 / 5 / $ -39 $ ,
4 / 1 / $ 3 $ ,4 / 2 / $ -13 $ ,4 / 3 / $ -29 $ ,4 / 4 / $ -45 $ ,4 / 5 / $ -61 $  }
    \latticepath{1}{(5,0)}{very thick}{
      -1/0/black, 0/-1/black,-1/0/black, 0/-1/black, 0/-1/black,-1/0/black, 0/-1/black, -1/0/black, -1/0/black}
  \end{tikzpicture}

  \caption{A lattice path in the array $A(8,11)$  }
  \label{PathInArray}
\end{figure}

\newcommand\MC[1]{\multicolumn{1}{c}{#1}}




To compute the average size of self-conjugate $(s,t)$-cores, we show that the size of a partition $\lambda$ can be expressed in terms of the entries in the array $A$ above the lattice path $P$ corresponding to $\lambda$ under the bijection $\Phi$.

\begin{lem}
For any lattice path $P$ in $\mathcal{P}(A)$, we have
\begin{displaymath}{\label{negtrick}}
|\Phi(P)|=  \frac{(s^2-1)(t^2-1)}{24}- \sum_{(i,j) \text{ is above } P}A_{i,j}.
\end{displaymath}
\end{lem}

\proof Clearly,  the size of a self-conjugate partition equals the sum of  the main diagonal hook lengths. By Theorem \ref{thmCorner}, we find that
\begin{eqnarray*}
|\Phi(P)|& = & \sum_{h\in MD(\Phi(P))} h\\[9pt]
&=& \sum_{(i,j) \text{ below } P  ,\,A_{i,j}>0} A_{i,j} - \sum_{(i,j) \text{ is above } P ,\, A_{i,j}<0} A_{i,j}\\ [9pt]
&=& \sum_{(i,j)\colon A_{i,j}>0} A_{i,j} - \sum_{(i,j) \text{ is above } P}A_{i,j}.
\end{eqnarray*}
To show that
\begin{equation}\label{231}
\sum_{(i,j)\colon A_{i,j}>0} A_{i,j}=\frac{(s^2-1)(t^2-1)}{24},
\end{equation}
let $Q$ be the lattice path along the left and upper borders of $A$.
Clearly, $M_A(Q)$ consists of positive entries of $A$. Let $\lambda=\Phi(Q)$.
By Theorem \ref{thmCorner},
the set of main diagonal hook length of $\lambda$ equals $M_A(Q)$.
Hence we obtain
\begin{equation}\label{232}
|\lambda| = \sum_{(i,j)\colon A_{i,j}>0} A_{i,j}.
\end{equation}

It remains to show that
\begin{equation}\label{233}
|\lambda| = \frac{(s^2-1)(t^2-1)}{24}.
\end{equation}

We claim that $\lambda$ is the largest $(s,t)$-core. Thus (\ref{233}) follows
from  the expression (\ref{largestSize}).
To prove this claim, we recall that Theorem \ref{largest} guarantees that there is an $(s,t)$-core with largest size, say $\mu$, that happens to be self-conjugate. We aim to show that $\mu=\lambda$. Let $l(\lambda)$ and $l(\mu)$ denote the lengths of
$\lambda$ and $\mu$ respectively. By Theorem \ref{thmCorner}, there is a lattice path $R \in \mathcal{P}(A)$ such that $\mu=\Phi(R)$. By Theorem \ref{contain}, we find that
\begin{equation}\label{mu_lambda2}
l(\mu)\geq l(\lambda)
\end{equation}
and
\begin{equation}\label{mu_lambda1}
\mu_i \geq \lambda_i
\end{equation}
for all $i$.
Combining (\ref{mu_lambda2}) and (\ref{mu_lambda1}), we obtain that
\begin{equation} \label{compareHook}
\mu_1 + l(\mu) -1 \geq \lambda_1 + l(\lambda) -1.
\end{equation}

The largest main diagonal hook length of $\lambda$ is $\lambda_1 + l(\lambda)-1$, that is,
\begin{equation}\label{maxLambda}
\lambda_1 + l(\lambda)-1 =\max MD(\lambda).
\end{equation}
Since $\lambda = \Phi(Q)$,  by Theorem \ref{thmCorner},  we have
\begin{equation}\label{MDLambda}
MD(\lambda)= M_A(Q) = \{ A_{i,j}| A_{i,j}>0, 1\leq i \le \fl{s/2},1\leq j \leq \fl{t/2} \}.
\end{equation}
Note that $A_{1,1}$ is the largest  in all positive entries in $A$. Thus, we deduce from (\ref{maxLambda}) and (\ref{MDLambda}) that
\begin{equation}\label{lam}
\lambda_1 + l(\lambda)-1=A_{1,1}.
\end{equation}

Since $\mu_1 + l(\mu)-1$  is the hook length of the cell in the upper-left corner of $\mu$,  by Theorem \ref{thmCorner}, $\mu_1 + l(\mu)-1$ belongs to $M_A(R)$.
By the definition of $M_A(R)$, there exists an entry $A_{i,j}$ of $M_A(R)$ such that
\begin{equation}\label{mu}
\mu_1 + l(\mu)-1 = |A_{i,j}|.
\end{equation}
We claim that \begin{equation} \label{11ij}
A_{1,1} \geq |A_{i,j}|,
 \end{equation}
for any entry $A_{i,j}$. Note that $A_{1,1}$ is the largest entry in $A$.
On the other hand, $A_{\fl{s/2},\fl{t/2}}$ is negative and is the smallest entry in $A$. It can be easily seen that  $|A_{\fl{s/2},\fl{t/2}}|< A_{1,1}$, since  \begin{equation*}\label{a11rules}
A_{1,1}+A_{\fl{s/2},\fl{t/2}} = st-s-t + st+s+t-2t\fl{s/2}-2s\fl{t/2}>0.
\end{equation*}
This proves the claim.

Combining \eqref{lam}, \eqref{mu} and \eqref{11ij}, we obtain that
\begin{equation} \label{compareHook2}
\lambda_1 + l(\lambda)-1 \geq \mu_1 + l(\mu)-1.
\end{equation}
From  (\ref{compareHook}) and (\ref{compareHook2}), we deduce that
\begin{equation}\label{eq1}
\lambda_1 + l(\lambda)-1 = \mu_1 + l(\mu)-1.
\end{equation}
By \eqref{lam} and \eqref{eq1}, we see that $A_{1,1}=\mu_1 + l(\mu)-1$, and hence
it is a main diagonal hook length of $\mu$. Thus $A_{1,1}$ lies in $MD(\mu)$. By Theorem \ref{thmCorner},  $A_{1,1}$ belongs to $M_A(R)$. Since $A_{1,1}>0$,
it is an entry of $A$ that is below the lattice path $R$. This implies that
$R$ is the unique lattice path of $A$ along the left and upper borders.
It follows that $Q=R$ and $\lambda=\mu$. So we conclude that $\lambda$ is the largest $(s,t)$-core. This completes the proof.

\qed

To prove the main result, we need some identities on the number of
lattice paths in a rectangular region. Let $m$ and $n$ be positive integers,
and $B_{mn}$ be an $m\times n$ diagram, that is, a diagram
of $m$ rows with each containing $n$ cells. The positions of the
cells of the first row are $(1,1), (1,2), \ldots,(1,n)$, and so on.
The set of lattice paths from the lower-left corner to the upper-right corner of
 $B_{mn}$ is denoted by $\mathcal{P}(B_{mn})$.
Let $f(i,j)$ be the number of lattice paths in $\mathcal{P}(B_{mn})$ that lie below the cell $(i,j)$, possibly touching the right or lower border of the cell $(i,j)$.

\begin{lem}\label{fij}
For positive integers $m,n$, we have
\begin{equation}\label{lem24}
        \sum_{1\leq i \leq m,\, 1\leq j \leq n} f(i,j)
    =   {{m+n}\choose m}  \frac{ m n }{2}.
\end{equation}
\end{lem}

\proof Given $1\leq i\leq m$ and $1\leq j\leq n$,  the number of lattice paths in $\mathcal{P}(B_{mn})$ below the cell $(i,j)$ equals the number of lattice paths above the cell $(m-i+1,n-j+1)$. Since each lattice path $P$ is either above the cell $(i,j)$ or below the cell $(i,j)$, we have
\begin{equation*}\label{24a}
f(i,j)+f(m-i+1,n-j+1)=|\mathcal{P}(B_{mn})|.
\end{equation*}
But the number of lattice paths in $\mathcal{P}(B_{mn})$ is ${{m+n}\choose m}$, we get
\begin{equation}\label{24b}
f(i,j)+f(m-i+1,n-j+1)={{m+n}\choose m}.
\end{equation}
Summing  \eqref{24b} over  $(i,j)$ gives $$2 \sum_{1\leq i \leq m,1\leq j \leq n} f(i,j)={{m+n}\choose m}mn.$$ This completes the proof.
\qed

\begin{lem}\label{ifij}
For positive integers $m$ and $n$, we have
\begin{equation}\label{l1}
    \sum_{1\leq i \leq m,\, 1\leq j \leq n} i   f(i,j)
=   {m+2 \choose 3}{m+n \choose {m+1}}
\end{equation}
and
\begin{equation}\label{l2}
    \sum_{1\leq i \leq m,\, 1\leq j \leq n} j  f(i,j)
=   {n+2 \choose 3}{m+n \choose {n+1}}.
\end{equation}
\end{lem}

\proof Let \[ G(m,n)=\sum_{1\leq i \leq m,\, 1\leq j \leq n} i  f(i,j).\] To prove  \eqref{l1}, we establish a  recurrence relation for $m,n\geq 2$,
\begin{equation}\label{lrec}
G(m,n)=G(m-1,n)+G(m,n-1)+{m+1 \choose 2}{m+n-1\choose m}.
\end{equation}
In doing so, let $T$ be the set of triples $(P,C_1,C_2)$, where $P$ is a path in $\mathcal{P}(B_{mn})$ , $C_1$ and $C_2$ are cells above $P$ and are in a same column with $C_2$ not lower than $C_1$. Notice that $C_1$ and $C_2$ are allowed to the
same cell.

We proceed to compute $|T|$ in two ways.
First, it is easily seen that $if(i,j)$ is the number of triples in $T$ with $C_1=(i,j)$. Hence we have for $m,n\geq 1$, $|T|=G(m,n)$.

Alternatively, $|T|$ can be computed as follows.

For a given lattice path $P$ in $\mathcal{P}(B_{mn})$, the cells above $P$ form a Ferrers diagram of a partition, denoted by $\mu$. Let $\mu'$ be the conjugate of $\mu$,
that is, there are $\mu'_j$ cells in the $j$-th column of the Ferrers diagram of $\mu$.

In the $j$-th column of the Ferrers diagram of $\mu$, there are ${{{\mu'_j}+1} \choose 2}$ ways to choose  $C_1$ and $C_2$ such that $C_2$ is not lower than $C_1$.
It follows that for given $P$, there are $\sum_{1\leq j\leq \mu_1} {\mu'_j +1 \choose 2}$ choices for $C_1$ and $C_2$.
Consequently, for $m,n\geq 1$,
\begin{equation}\label{l32}
|T|=\sum_{\mu\colon1\leq  \mu_1 \leq n,\, 1\leq \mu'_1 \leq m} \;\sum_{1\leq j\leq \mu_1} {\mu'_j +1 \choose 2}.
\end{equation}
Hence, for $m,n\geq 1$,
\begin{equation}\label{rr1}
G(m,n)=\sum_{\mu\colon1\leq  \mu_1 \leq n,\, 1\leq \mu'_1 \leq m} \;\sum_{1\leq j\leq \mu_1} {\mu'_j +1 \choose 2}.
\end{equation}

%
For $m,n\geq 2$, the right hand side of \eqref{rr1} equals
\begin{equation}\label{rr2}
\sum_{\mu\colon1\leq \mu_1 \leq n,\, \mu'_1= m} \; \sum_{1\leq j\leq \mu_1} {\mu'_j +1 \choose 2}
 + \sum_{\mu\colon 1\leq \mu_1 \leq n,\, 1\leq \mu'_1 \leq m-1} \;\sum_{1\leq j\leq \mu_1} {\mu'_j +1 \choose 2}.
\end{equation}
It is evident from \eqref{rr1} that the second double sum in \eqref{rr2} can be expressed by $G(m-1,n)$. The first double sum in \eqref{rr2} can be rewritten as
\begin{equation}\label{rr3}
\sum_{ \mu\colon 1\leq \mu_1 \leq n,\, \mu'_1 = m}\; \sum_{2\leq j\leq \mu_1} {{\mu}'_j +1 \choose 2}+\sum_{ \mu\colon 1\leq \mu_1 \leq n,\, \mu'_1 = m}\; {{m+1} \choose 2}.
\end{equation}
Clearly, the number of partitions $\mu$ with $1\leq \mu_1 \leq n$ and $\mu'_1 = m$
 equals the number of lattice paths from the lower-left corner
 to the upper-right corner in $B_{m,n-1}$, which is ${{m+n-1}\choose m}$. Hence
the second sum in \eqref{rr3} simplifies to
\begin{equation}\label{rr4}
{{m+1} \choose 2}{{m+n-1}\choose m}.\end{equation}

To compute the double sum in \eqref{rr3}, let $\tilde{\mu}$ denote the partition obtained from $\mu$ by deleting the first column of the Ferrers diagram of $\mu$.
So we see that
\begin{eqnarray} \label{Gmn-1}
\sum_{ \mu\colon 1\leq \mu_1 \leq n,\, \mu'_1 = m} \; \sum_{2\leq j\leq \mu_1} {{\mu}'_j +1 \choose 2}
&=&
\sum_{ \tilde{\mu}\colon 0 \leq \tilde{\mu}_1 \leq n-1,\, \tilde{\mu}'_1 \leq m}\;  \sum_{1\leq j\leq \tilde{\mu}_1} {\tilde{\mu}'_j +1 \choose 2}.
\end{eqnarray}
From \eqref{rr1} it can be seen that the right hand side of \eqref{Gmn-1} equals $G(m ,n-1)$.
Combining \eqref{rr1}--\eqref{Gmn-1}, we arrive at
the recurrence relation \eqref{lrec}.

For $m,n\geq 1$, let
$$F(m,n)= {m+2 \choose 3} {m+n\choose m+1}.$$
To prove that $G(m,n)=F(m,n)$ for $m,n \geq 1$,
 it is sufficient to check that
  $F(m,n)$ satisfies the same recurrence relation \eqref{lrec} and
  the same initial conditions.
Clearly, $F(1, n)=G(1, n)$ and $F(m,1)=G(m,1)$ for $m,n\geq 1$.
Moreover, it is easily checked that the
recurrence relation \eqref{lrec} holds for $F(m,n)$ as well.
This proves identity \eqref{l1}.
Relation  \eqref{l2} can be viewed as a restatement of (\ref{l1}).
This completes the proof.
 \qed

Now we are ready to prove the conjecture of Armstrong, Hanusa and Jones
on the average size of a self-conjugate $(s,t)$-core.

\proof Let $SC(s,t)$ denote the set of self-conjugate $(s,t)$-cores.  We
aim to show that
\begin{equation}\label{mm}
\sum_{\lambda \in SC(s,t)} |\lambda| = \frac{(s+t+1)(s-1)(t-1)}{24} {{\fs+\ft}\choose \fs}.
\end{equation}
By Theorem \ref{thmCorner}, we find that
\begin{equation} \label{avesize}
\sum_{\lambda \in SC(s,t)} |\lambda| = 
\sum_{P\in\mathcal{P}(A)} |\Phi(P)|.
\end{equation}
Using Lemma \ref{negtrick}, we obtain that
\begin{equation}\label{as1}
\sum_{P\in\mathcal{P}(A)} |\Phi(P)|=\frac{(s^2-1)(t^2-1)}{24} {{\fs+\ft}\choose \ft}-\sum_{P\in\mathcal{P}(A)} \;\sum_{(i,j) \text{ is above } P}A_{i,j}.
\end{equation}
Combining \eqref{avesize} and \eqref{as1}, we see that
\begin{equation}\label{eq:minus}
\sum_{\lambda \in SC(s,t)} |\lambda|=\frac{(s^2-1)(t^2-1)}{24} {{\fs+\ft}\choose \ft}-\sum_{P\in\mathcal{P}(A)}\; \sum_{(i,j) \text{ is above } P}A_{i,j}.
\end{equation}
By the definition \eqref{aij} of the array $A$, we deduce that
\begin{eqnarray}\label{avesize1}
 \sum_{P\in\mathcal{P}(A)}\; \sum_{(i,j) \text{ is above } P}A_{i,j}
&=&\sum_{P\in\mathcal{P}(A)}\; \sum_{(i,j) \text{ is above } P} (st+s+t-2sj-2ti) \nonumber \\[9pt]
&=&
(st+s+t)\sum_{1\leq i\leq \fs,\, 1\leq j\leq \ft} f(i,j)  -2s\sum_{1\leq i\leq \fs,\, 1\leq j\leq \ft} j f(i,j)\nonumber\\[9pt]
& &\quad -2t\sum_{1\leq i\leq \fs,\, 1\leq j\leq \ft}i f(i,j).
\end{eqnarray}
Applying Lemma \ref{fij} and  Lemma \ref{ifij}  to (\ref{avesize1}) with   $m=\fs$ and $n=\ft$, we get
\begin{eqnarray}
\sum_{P\in\mathcal{P}(A)}\; \sum_{(i,j) \text{ is above } P}A_{i,j}&=&
(st+s+t)  {{\fs+\ft}\choose \fs}  \frac{ \fs \ft}{2}
-2s  {\ft+2 \choose 3}{\fs+\ft \choose {\fs-1}}\nonumber\\
& &\quad -2t  {\fs+2 \choose 3}{\fs+\ft \choose {\ft-1}} .\label{avesize3}
\end{eqnarray}
We claim that
\begin{eqnarray*}
\frac{(s^2-1)(t^2-1)}{24}   {{\fs+\ft}\choose \ft} & = &\frac{(s+t+1)(s-1)(t-1)}{24} {{\fs+\ft}\choose \ft} \\& &\quad +(st+s+t){\frac{\fs \ft}{2}}{{\fs+\ft}\choose \ft} \\& &\quad -2t {{\fs+2}\choose 3}{{\fs+\ft}\choose \ft-1} -2s {{\ft+2}\choose 3}{{\fs+\ft}\choose \fs-1},
\end{eqnarray*}
which simplifies to
\begin{eqnarray*}
\frac{st(s-1)(t-1)}{24}=(st+s+t){\frac{\fs \ft} 2}-{\frac t 3}(\fs+2)\fs\ft-{\frac s 3}(\ft+2)\fs\ft.
\end{eqnarray*}
When $s$ and $t$ are coprime, at least one of $s$ and $t$ is odd.
Without loss of generality, we may assume that $s$ is odd. In this case,
it is easily checked that above relation is true. Thus the claim holds. Combining \eqref{eq:minus} and \eqref{avesize3},
we arrive at \eqref{mm}, and hence the proof is complete.
\qed



%


%


\vspace{.2cm} \noindent{\bf Acknowledgments.}
This work was supported by the 973 Project, the PCSIRT Project of the Ministry of Education   and the National Science Foundation of China.


\end{document}